\documentclass{article}
\usepackage{latexsym}
\usepackage{graphics}
\usepackage{amsfonts}
\newtheorem{thm}{Theorem}
\newtheorem{conj}{Conjecture}
\newtheorem{prop}{Proposition}

\newtheorem{defn}{Definition}
\newtheorem{cor}{Corollary}

\begin{document}
\begin{center}
{\bf A FEW WEIGHT SYSTEMS ARISING FROM INTERSECTION GRAPHS}\\
\vspace{.2in}
{\footnotesize BLAKE MELLOR}\\
{\footnotesize Loyola Marymount University}\\
{\footnotesize 7900 Loyola Boulevard}\\
{\footnotesize Los Angeles, CA  90045-8130}\\
{\footnotesize\it  bmellor@lmu.edu}\\
\vspace{0.5in}
{\footnotesize ABSTRACT}\\
{\ }\\
\parbox{4.5in}{\footnotesize \ \ \ \ \ We show that the adjacency matrices of the
intersection graphs of chord diagrams satisfy the 2-term relations
of Bar-Natan and Garoufalides \cite{bg}, and hence give rise to
weight systems.  Among these weight systems are those associated
with the Conway and HOMFLYPT polynomials.  We extend these ideas
to looking at a space of {\it marked} chord diagrams modulo an
extended set of 2-term relations, define a set of generators for
this space, and again derive weight systems from the adjacency
matrices of the (marked) intersection graphs. Among these weight
systems are those associated with the Kauffman polynomial.}\\
\vspace{1in}
\end{center}
\tableofcontents

\section{Introduction} \label{S:intro}

Finite type invariants have received a lot of attention over the
past decade.  One reason for this is that they provide a common
framework for many of the most powerful knot invariants, such as
the Conway, Jones, HOMFLYPT and Kauffman invariants.  The
framework also allows us to study these invariants using
elementary combinatorics, by looking at associated functionals
(called {\it weight systems}) on spaces of {\it chord diagrams}.
This provides a new ways of describing the invariants.

The modest goal of this paper is to define a few weight systems in
terms of the adjacency matrix of the intersection graph of the
chord diagrams, and to show that among these weight systems are
those associated with the Conway, HOMFLYPT and Kauffman
polynomials, in both their framed and unframed incarnations.  This
gives us new formulas for the weight systems associated to these
important knot invariants.  We build on ideas of Bar-Natan and
Garoufalides \cite{bg}, who first found the formula we give for
the Conway polynomial.

In section~\ref{S:prelim} we will review the necessary background
for the paper, including finite type invariants, the 2-term
relations introduced by Bar-Natan and Garoufalides, intersection
graphs of chord diagrams and Lando's graph bialgebra.  In
section~\ref{S:adj} we will study the adjacency matrix of the
intersection graph, and show that the weight systems associated
with the Conway and HOMFLYPT polynomials can be defined in terms
of the determinant and rank of this matrix.  In
section~\ref{S:marked} we look at {\it marked} chord diagrams and
define an extended set of 2-term relations on these diagrams.  We
give an explicit set of generators for the space of marked chord
diagrams modulo these relations.  Finally, we show that the weight
system associated with the Kauffman polynomial can be defined in
terms of the rank of the adjacency matrix of marked chord
diagrams.

{\it Remark:}  The result for the Conway polynomial
(Theorem~\ref{T:conway}) has been previously proved by Bar-Natan
and Garoufalidis \cite{bg}, but is included here for completeness
and to place it in the context of Lando's bialgebra. After
distributing the first version of this paper \cite{me2}, the
author discovered that the adjacency matrix of an intersection
graph has also been studied by Soboleva \cite{so}, who has also
proven Theorem~\ref{T:homfly}, and a weaker version of
Theorem~\ref{T:kauff}. The intersection graphs we study are also
related to the {\it trip matrix} of a knot, studied by Zulli
\cite{zu}.

\section{Preliminaries} \label{S:prelim}

\subsection{Finite Type Invariants} \label{SS:finitetype}

In 1990, V.A. Vassiliev introduced the idea of {\it Vassiliev} or
{\it finite type} knot invariants, by looking at certain groups
associated with the cohomology of the space of knots.  Shortly
thereafter, Birman and Lin~\cite{bl} gave a combinatorial
description of finite type invariants.  We will give a brief
overview of this combinatorial theory.
  For more details, see Bar-Natan~\cite{bn}.

A {\it knot} is an embedding of the circle $S^1$ into the 3-sphere
$S^3$.  A {\it knot invariant} is a map from these embeddings to
some set which is invariant under isotopy of the embedding.  We
will also consider invariants of {\it regular} isotopy, where the
isotopy preserves the {\it framing} of the knot (i.e. a chosen
section of the normal bundle of the knot in $S^3$).  We first note
that we can extend any knot invariant to an invariant of {\it
singular} knots, where a singular knot is an immersion of $S^1$ in
3-space which is an embedding except for a finite number of
isolated double points. Given a knot invariant $v$, we extend it
via the relation:

$$\includegraphics{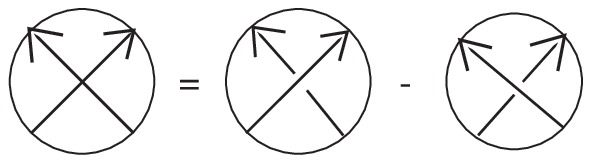}$$

An invariant $v$ of singular knots is then said to be of {\it
finite type}, specifically of {\it type n}, if $v$ is zero on any
knot with more than $n$ double points (where $n$ is a finite
nonnegative integer). The smallest such $n$ is called the {\it
order} of $v$. We denote by $V_n$ the vector space over ${\mathbb
C}$ generated by (framing-independent) finite type invariants of
type $n$ (i.e., whose order is $\leq$ $n$).  We can completely
understand the space of finite type invariants by understanding
all of the vector spaces $V_n/V_{n-1}$.  An element of this vector
space is completely determined by its behavior on knots with
exactly $n$ singular points.  In addition, since such an element
is zero on knots with more than $n$ singular points, any other
(non-singular) crossing of the knot can be changed without
affecting the value of the invariant.  This means that elements of
$V_n/V_{n-1}$ can be viewed as functionals on the space of {\it
chord diagrams}:
\begin{defn}
A {\bf chord diagram of degree n} is an oriented circle, together with $n$
chords of the circles, such that all of the $2n$ endpoints of the chords are
distinct.  The circle represents a knot, the endpoints of a chord represent
2 points identified by the immersion of this knot into 3-space.  The diagram
is determined by the order of the $2n$ endpoints.
\end{defn}

Functionals on the space of chord diagrams which are derived from
finite type knot invariants will satisfy certain relations.  This
leads us to the definition of a {\it weight system}:
\begin{defn}
A {\bf weight system of degree n} is a linear functional $W$ on
the space of chord diagrams of degree $n$ (with values in an
associative commutative ring ${\bf K}$ with unity) which satisfies
the 1-term and 4-term relations, shown in Figure~\ref{F:4-term}.
    \begin{figure}
    (1-term relation) $$\includegraphics{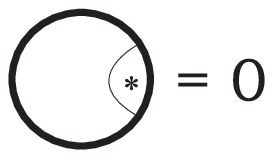}$$
    (4-term relation) $$\includegraphics{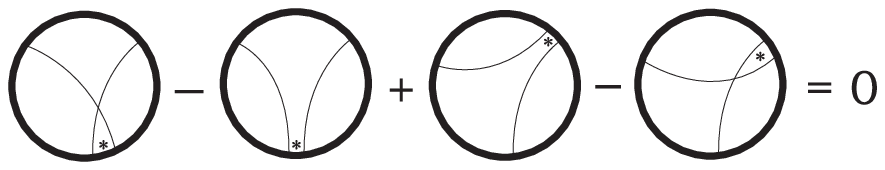}$$
    \caption{The 1-term and 4-term relations.  No other chords have endpoints
    on the arcs labeled with a *.  In the 4-term relations, all other chords of the four
    diagrams are the same.} \label{F:4-term}
    \end{figure}
\end{defn}

It can be shown (see \cite{bl,bn,st}) that the space of weight
systems of degree $n$ is isomorphic to $V_n/V_{n-1}$.  For
convenience, we take the dual approach, and simply study the space
of chord diagrams of degree $n$ modulo the 1-term and 4-term
relations. The 1-term relation is occasionally referred to as the
"framing-independence" relation, because it arises from the
framing-independence of the invariants in $V_n$ (essentially, from
the first Reidemeister move). Since most of the interesting
structure of the vector spaces arises from the 4-term relation, it
is common to look at the more general setting of invariants of
regular isotopy, and consider the vector space $A_n$ of chord
diagrams of degree $n$ modulo the 4-term relation alone.  We will
call the space $W_n$ of linear functionals on $A_n$ the space of
{\it regular weight systems} of degree $n$.  We will let
$\widehat{A_n}$ denote the vector space of chord diagrams modulo
both the 1-term and 4-term relations, and $\widehat{W_n}$ denote
the space of functionals on $\widehat{A_n}$, the unframed weight
systems.

It is useful to combine all of these spaces into a graded module
$A = \bigoplus_{n\geq 1}A_n$ via direct sum.  We can give this
module a bialgebra (or Hopf algebra) structure by defining an
appropriate product and co-product:
\begin{itemize}
    \item  We define the product $D_1 \cdot D_2$ of two chord diagrams
$D_1$ and $D_2$ as their connect sum.  This is well-defined modulo
the 4-term relation (see \cite{bn}).
$$\includegraphics{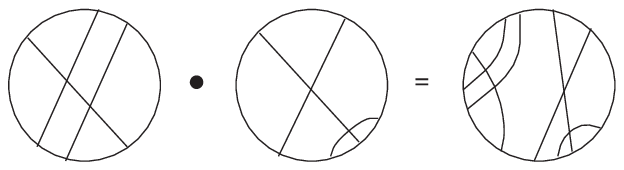}$$
    \item  We define the co-product $\Delta(D)$ of a chord diagram $D$ as
follows:
$$\Delta(D) = {\sum_J D_J' \otimes D_J''}$$
where $J$ is a subset of the set of chords of $D$, $D_J'$ is $D$
with all the chords in $J$ removed, and $D_J''$ is $D$ with all
the chords {\it not} in $J$ removed.
$$\includegraphics{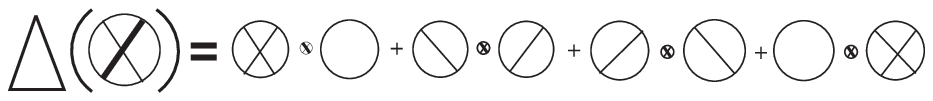}$$
\end{itemize}
It is easy to check the compatibility condition $\Delta(D_1\cdot
D_2) = \Delta (D_1)\cdot\Delta(D_2)$.

There is a natural {\it deframing} map $\phi: A \otimes A
\rightarrow A$, defined by:
$$\phi(D_1 \otimes D_2) = (-\Theta)^{deg(D_1)} \cdot D_2$$
Here $\Theta$ represents the chord diagram consisting of a single
chord.  This map gives a canonical projection $\widehat{}:W_n
\rightarrow \widehat{W_n}$, defined by
$\widehat{W}(D)=W(\phi(\Delta(D)))$ (see \cite{bn}, Exercise
3.16).

\subsection{2-Term relations} \label{SS:2term}

\begin{figure}
    $$\includegraphics{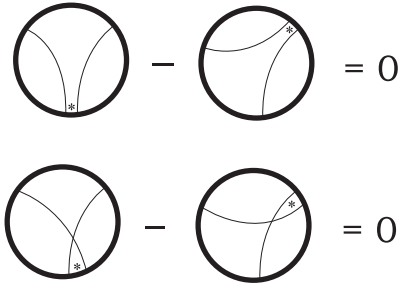}$$
    \caption{The 2-term relations} \label{F:2-term}
    \end{figure}

Of course, any particular weight system will satisfy relations in
addition to the 1-term and 4-term relations, and it can be useful
to look at weight systems which lie in the subspaces determined by
these additional relations.  In particular, Bar-Natan and
Garoufalides \cite{bg} noted that the weight system associated
with the Conway polynomial satisfies the {\it 2-term} relations in
Figure~\ref{F:2-term}.
Clearly, these relations imply that the weight system satisfies
the 4-term relation as well.  As a result, the product and
coproduct of section \ref{SS:finitetype} are still well-defined.
So we can give the vector space of chord diagrams modulo the
2-term relations the structure of a bialgebra.  We will denote
this bialgebra (and the underlying vector space) by $B$.  There is
a natural projection from $A$ to $B$.

\begin{figure}
    $$\includegraphics{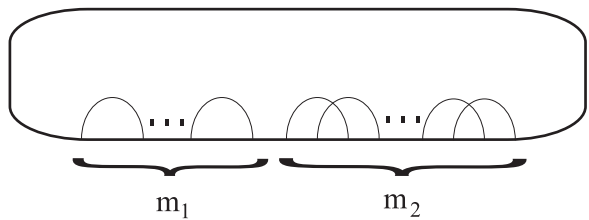}$$
    \caption{Example of an ($m_1, m_2$)-caravan} \label{F:caravan}
    \end{figure}

Bar-Natan and Garoufalides also showed that $B$ is generated (as a
vector space) by {\it ($m_1,m_2$)-caravans} of $m_1$ "one-humped
camels" (isolated chords which intersect no other chords) and
$m_2$ "two-humped camels" (pairs of chords which intersect each
other, but no other chords). An example of such a caravan is shown
in Figure~\ref{F:caravan}.

\subsection{Intersection Graphs} \label{SS:graphs}

\begin{defn}
Given a chord diagram $D$, we define its {\bf intersection graph} $\Gamma(D)$ as the graph
such that:
\begin{itemize}
    \item $\Gamma(D)$ has a vertex for each chord of $D$.
    \item Two vertices of $\Gamma(D)$ are connected by an edge if and only if the corresponding
chords in $D$ intersect, i.e. their endpoints on the bounding circle alternate.
\end{itemize}
\end{defn}
For example:
$$\includegraphics{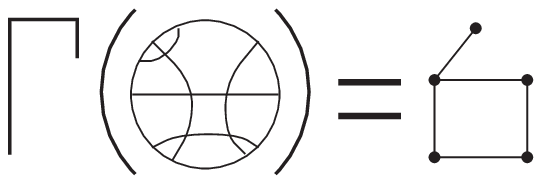}$$
Note that these graphs are simple (i.e. every edge has two
distinct endpoints, and there is at most one edge connecting any
two vertices). These graphs are also known as {\it circle graphs},
and have been studied extensively by graph theorists.  A
combinatorial classification of circle graphs has been given by
Bouchet \cite{bo}.

A circle graph can be the intersection graph for more than one chord diagram.  For example, there
are three different chord diagrams with the following intersection graph:
$$\includegraphics{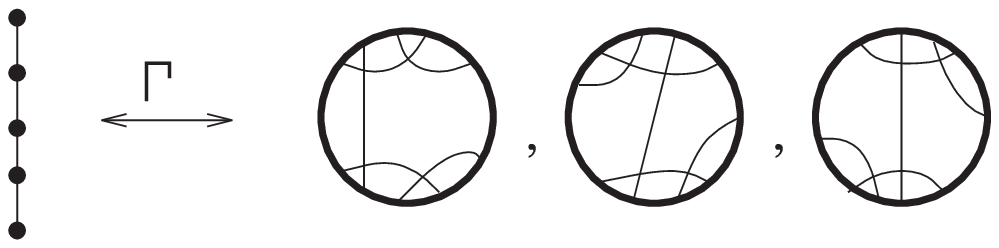}$$

However, these chord diagrams are all equivalent modulo the 4-term relation.  Chmutov,
Duzhin and Lando~\cite{cdl} conjectured that intersection graphs actually determine
the chord diagram, up to the 4-term relation.  In other words, they proposed:
\begin{conj} \label{C:IGC}
If $D_1$ and $D_2$ are two chord diagrams with the same intersection graph, i.e. $\Gamma(D_1)$
= $\Gamma(D_2)$, then for any weight system $W$, $W(D_1) = W(D_2)$.
\end{conj}
This Intersection Graph Conjecture is now known to be false in general.  Morton
and Cromwell~\cite{mc} found a finite type invariant of type 11 which can distinguish some mutant knots,
and Le~\cite{le} and Chmutov and Duzhin~\cite{cd} have shown that mutant knots cannot be distinguished by
intersection graphs.  However, the conjecture is true in many special cases, and the exact extent to which
it fails is still unknown, and potentially very interesting.

The conjecture is known to hold in the following cases \cite{cdl}:
\begin{itemize}
    \item  For chord diagrams with 8 or fewer chords (checked via computer calculations);
    \item  For the weight systems coming from the defining representations of Lie
algebras {\it gl(N)} or {\it so(N)} as constructed by Bar-Natan in \cite{bn};
    \item  When $\Gamma(D_1) = \Gamma(D_2)$ is a tree (or, more generally, a linear
combination of forests);
    \item  When $\Gamma(D_1) = \Gamma(D_2)$ has a single loop (see \cite{me}).
\end{itemize}

The second item above includes the weight systems arising from the
Conway, HOMFLYPT and Kauffman polynomials.  A main goal of this
paper is to find explicit formulas for these weight systems in
terms of intersection graphs.

\subsection{Lando's graph bialgebra} \label{SS:Lando}

Lando \cite{la} has given more structure to the questions
surrounding intersection graphs by extending the map $\Gamma$ to a
homomorphism between the bialgebra $A$ of chord diagrams and a
particular bialgebra of graphs. Lando's bialgebra of graphs is
constructed by defining an analogue of the 4-term relation for
graphs, as follows:

\begin{defn} \cite{la}
Consider the graded vector space (over ${\mathbb C}$) of formal
linear combinations of graphs, graded by the number of vertices in
the graphs.  For any graph G and vertices A and B in V(G) we
impose on the vector space the relation:
$$G-G'_{AB}-\widetilde{G}_{AB}+\widetilde{G}'_{AB} = 0$$
where $G'_{AB}$ is the result of complementing the edge AB in G
(i.e. adding or removing it), $\widetilde{G}_{AB}$ is the result
of complementing the edge AC for every vertex C in V(G) which is
adjacent to B and $\widetilde{G}'_{AB}$ is the result of
complementing the edge AB in $\widetilde{G}_{AB}$.  Here is an
example of such a relation:
$$\includegraphics{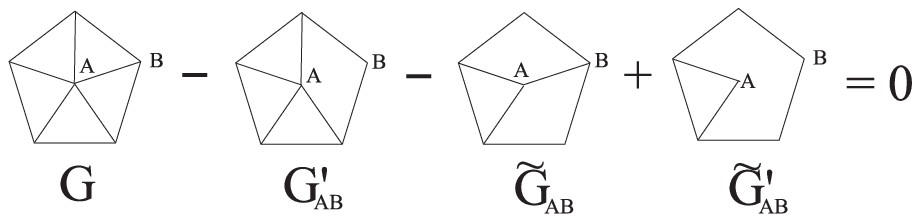}$$
The bialgebra $F$ is defined as this graded vector space, together
with a product and a coproduct.  The product is simply disjoint
union of graphs.  The coproduct is a map $\mu: F \rightarrow F
\otimes F$, defined as follows.  For any graph G, and subset $J
\subseteq V(G)$ of its vertices, let $G_J$ denote the subgraph
induced by $J$.  Then:
$$\mu(G) = \sum_{J \subseteq V(G)}{G_J \otimes G_{V(G)\backslash J}}$$
An example is shown below:
$$\includegraphics{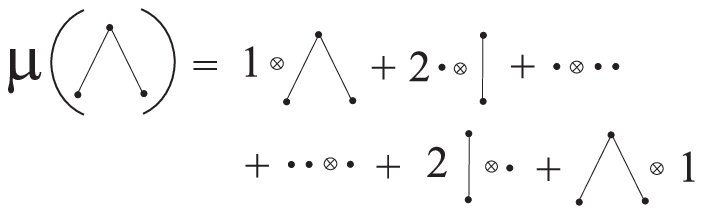}$$
\end{defn}

It is now easy to show that $\Gamma$ extends to a bialgebra
homomorphism from $A$ to $F$ (see \cite{la}).

We can easily extend Lando's results to include the 1-term
relation and framing-independent invariants.  We define the
algebra $\widehat{F}$ to be simply $F$ modulo graphs with isolated
vertices (these correspond to the isolated chords of the 1-term
relation for chord diagrams).  It is then trivial to show that
$\Gamma$ extends to a bialgebra homomorphism from $\widehat{A}$ to
$\widehat{F}$.

A {\it regular graph weight system} is a linear functional
$\gamma: F \rightarrow {\mathbb C}$ (Lando called these
functionals {\it 4-invariants}) . Then, given any regular graph
weight system $\gamma$, $\gamma\circ\Gamma: A \rightarrow {\mathbb
C}$ is a regular weight system. Similarly, if we define a {\it
graph weight system} to be a linear functional of $\widehat{F}$,
then for any graph weight system $\alpha$, $\alpha\circ\Gamma$
will be a weight system.

Just as for chord diagrams, there is a natural {\it deframing} map
$\phi: F \otimes F \rightarrow F$, defined by:
$$\phi(G_1 \otimes G_2) = (-\bullet)^{deg(G_1)} \cdot G_2$$
Here $\bullet$ represents the trivial graph consisting of a single
vertex and no edges. This map gives a canonical projection
$\widehat{}:F^* \rightarrow \widehat{F}^*$, defined by
$\widehat{\gamma}(G)=\gamma(\phi(\mu(G)))$.

\section{The Adjacency Matrix of an Intersection Graph} \label{S:adj}

In this section we will show that the determinant and rank of the
adjacency matrix of a graph (over ${\bf Z}_2$) are regular graph
weight systems, and that the determinant is, in addition, a graph
weight system.  We will do this by showing that the isomorphism
class of the adjacency matrix (as a symmetric bilinear form over
${\bf Z}_2$) satisfies 2-term relations analagous to those in
section \ref{SS:2term}.  We will then show that these weight
systems are essentially the same as those associated with the
Conway and HOMFLYPT polynomials.

\subsection{Graph Weight Systems from the Adjacency Matrix}
\label{SS:adjweight}

We begin by recalling the definition of the adjacency matrix of a
graph.

\begin{defn} \label{D:adj}
Given a graph G with n vertices, labeled $\{v_1,...,v_n\}$, the {\bf adjacency matrix of G}, or adj(G),
is the symmetric $n \times n$ matrix defined by:
$$adj(G)_{ij} = \left\{{\matrix{1\ if\ v_i\ and\ v_j\ are\ connected\ by\ an\ edge\ in\ G \cr 0\
otherwise}}\right.$$
In the case of a simple graph, the diagonal entries of the matrix will all be 0.
\end{defn}

This matrix can be viewed as a symmetric bilinear form over ${\bf
Z}_2$.  If we permute the labels on the vertices of $G$, we change
the matrix $adj(G)$ by the corresponding permutations of the rows
and columns. But this does not change the isomorphism class of the
form (see \cite{mh}).  So, as an isomorphism class of symmetric
bilinear forms, the adjacency matrix of an unlabeled graph is
well-defined. From Milnor and Husemoller \cite{mh}, we know that
the determinant and rank of the matrix are invariants of the
isomorphism class of the form, and hence are well-defined
invariants of the graph.  This leads us to define the following
functions on graphs:

\begin{defn}
Given a graph G, we define the {\bf determinant of G} and the {\bf
rank of G} as follows:
$$det(G) = det(adj(G)) \in {\bf Z}_2$$
$$rank(G) = rank(adj(G))$$
\end{defn}

We extend these functions linearly to get {\bf Z}-valued
functionals on the space of graphs.  We will also call these
extensions the determinant and rank.  We will see that the
determinant gives a {\bf Z}-valued graph weight system, and the
rank gives a {\bf Z}-valued regular graph weight system (the rank
does not satisfy the 1-term relation).  We first show that both
functionals are regular graph weight systems.  To do this, we will
show that they satisfy {\it 2-term} relations, analogous to those
in section \ref{SS:2term}, defined as follows. Consider graphs $G,
G'_{AB}, \widetilde{G}_{AB}, \widetilde{G}'_{AB}$ as in section
\ref{SS:Lando}.  Then the 2-term relations are:
$$G-\widetilde{G}_{AB}=0$$
$$G'_{AB}-\widetilde{G}'_{AB}=0$$
It is clear that any functional which satisfies these 2-term
relations will also satisfy the 4-term relation.  So the vector
space $E$ of graphs modulo the 2-term relations can be given the
structure of a bialgebra, using the same product and coproduct as
for $F$.  There is a natural projection from $F$ to $E$. Moreover,
the pullback by $\Gamma$ of any functional on $E$ will be a
functional on $B$ (defined in section \ref{SS:2term}).

\begin{thm} \label{T:weight}
The isomorphism class of the adjacency matrix of a graph satisfy
the 2-term relations above.
\end{thm}
{\sc Proof:}  Consider two vertices $A$ and $B$, giving rise to
the four graphs $G, G'_{AB}, \widetilde{G}_{AB},
\widetilde{G}'_{AB}$. We want to show that $adj(G) \cong
adj(\widetilde{G}_{AB})$ and $adj(G'_{AB}) \cong
adj(\widetilde{G}'_{AB})$. The easiest way to do this is simply to
write down the matrices explicitly.  The vertices of $G$ other
than $A$ and $B$ can be partitioned into four sets $S_{AB}, S_A,
S_B, {\rm \ and\ } S_0$, where $S_{AB}$ contains those vertices
adjacent to both $A$ and $B$ in $G$, $S_A$ contains those vertices
adjacent to $A$ but not $B$ in $G$, $S_B$ contains those vertices
adjacent to $B$ but not $A$ in $G$, and $S_0$ contains those
vertices adjacent to neither $A$ nor $B$ in $G$.

The adjacency matrices for the four graphs, with respect to the
basis $\{A, B, S_{AB}, S_A, S_B, S_0\}$, are shown below. We
assume that $A$ and $B$ are connected by an edge in $G$ (if not,
simply interchange $G$ and $G'$). Here $I$ and $O$ represent a row
or column of 1's and 0's respectively:
$$adj(G) = \left[{\matrix{0 & 1 & I & I & O & O \cr
    1 & 0 & I & O & I & O \cr
    I & I & * & * & * & * \cr
    I & O & * & * & * & * \cr
    O & I & * & * & * & * \cr
    O & O & * & * & * & *}}\right] \cong \left[{\matrix{0 & 1 & O & I & I & O \cr
    1 & 0 & I & O & I & O \cr
    O & I & * & * & * & * \cr
    I & O & * & * & * & * \cr
    I & I & * & * & * & * \cr
    O & O & * & * & * & *}}\right] = adj(\widetilde{G}_{AB})$$
$$adj(G'_{AB}) = \left[{\matrix{0 & 0 & I & I & O & O \cr
    0 & 0 & I & O & I & O \cr
    I & I & * & * & * & * \cr
    I & O & * & * & * & * \cr
    O & I & * & * & * & * \cr
    O & O & * & * & * & *}}\right] \cong \left[{\matrix{0 & 0 & O & I & I & O \cr
    0 & 0 & I & O & I & O \cr
    O & I & * & * & * & * \cr
    I & O & * & * & * & * \cr
    I & I & * & * & * & * \cr
    O & O & * & * & * & *}}\right] = adj(\widetilde{G}'_{AB})$$
The isomorphisms are just the result of adding the second row (and
column) of the matrix on the left to its first row (and column),
modulo 2. So the isomorphism classes of the adjacency matrices of
the graphs satisfy the 2-term relations.  $\Box$

\begin{cor} \label{C:rankdet}
The (linear extensions of) the rank and determinant of a graph are
regular graph weight systems.
\end{cor}

\begin{thm} \label{T:det5}
The (linear extension of) the determinant of a graph is a graph
weight system.
\end{thm}
{\sc Proof:}  We need to show that the determinant of a graph
satisfies the 1-term relation - i.e. that it is trivial on graphs
with isolated vertices.  Let $G$ be a graph with an isolated
vertex, $v$, and let $G* = G-\{v\}$.  Then the adjacency matrix
for $G$ can be represented:
$$adj(G) = \left[{\matrix{0 & 0 \cr 0 & adj(G*)}}\right]$$
Since there is a row (and column) of 0's, $det(G) = 0$, so the
determinant satisfies the 1-term relation.  $\Box$

As we mentioned earlier, the rank of a graph does {\it not}
satisfy the 1-term relation, so it is not a graph weight system.
However, we can use the canonical projection from section
\ref{SS:Lando} to construct a graph weight system from the rank.
In fact, we will construct a {\it polynomial} graph weight system
by beginning with the invariant $R(G)(x) = x^{rank(G)}$, whose
linear extension is clearly also a regular graph weight system. To
apply our projection, it suffices to note that $rank(G_1 \cdot
G_2) = rank(G_1) + rank(G_2)$, and so $rank(\bullet^{deg(G_1)}
\cdot G_2) = deg(G_1)rank(\bullet) + rank(G_2) = rank(G_2)$.

\begin{thm} \label{T:r(d)}
Given a graph G, we define a polynomial $\widehat{R}(G)(x)$ as
follows. Here J is a subset of the vertices of G, $|J|$ is the
size of J, n is the total number of vertices in G, and $G_J$ is
the subgraph of G induced by J:
$$\widehat{R}(G)(x) = \sum_J{(-1)^{n-|J|}x^{rank(G_J)}}$$
This polynomial is the canonical projection of R(G), and so its
linear extension to a {\bf Z}[x]-valued functional on the space of
graphs is a graph weight system.
\end{thm}

\subsection{The Conway and HOMFLYPT weight systems}
\label{SS:Conway}

\begin{figure}
    $$\includegraphics{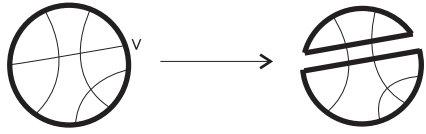}$$
    \caption{Surgery on a chord $v$} \label{F:surgery}
    \end{figure}

The Conway polynomial $\Delta$ of a link is a power series
$\Delta(L) = \sum_{n \ge 0}{a_n(L)z^n}$.  It can be computed via
the skein relation (where $L_+, L_-, L_0$ are as in
Figure~\ref{F:skein}):
    \begin{figure}
    $$\includegraphics{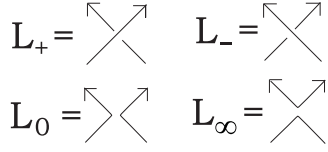}$$
    \caption{Diagrams of the skein relation.} \label{F:skein}
    \end{figure}
$$\Delta(L_+)-\Delta(L_-) = z\Delta(L_0)$$
$$\Delta(unlink\ of\ k\ components) = \left\{{\matrix{1\ if\ k=1 \cr 0\ if\ k>1}}\right.$$
The coefficient $a_n$ is a finite type invariant of type $n$ (see
\cite{bl}, \cite{bn}), and therefore defines a weight system $b_n$
of degree $n$.  The collection of all these weight systems is
called the Conway weight system, denoted $C$.  Consider a chord
diagram $D$, together with a chord $v$.  Let $D_v$ be the result
of {\it surgery on v}, i.e. replacing $v$ by an untwisted band,
and then removing the interior of the band and the intervals where
it is attached to $D$, as shown in Figure~\ref{F:surgery}
(so $D_v$ may have multiple boundary circles).  The skein
relations for the Conway polynomial give rise to the following
relations for $C$:
$$C(D) = C(D_v)$$
$$C(unlink\ of\ k\ components) = \left\{{\matrix{1\ if\ k=1 \cr 0\ if\ k>1}}\right.$$

It is easy to show (see \cite{bg}) that this weight system
satisfies the 2-term relations of section \ref{SS:2term}.  Simply
surger the two chords; the 2-term relation then says just that one
band can be "slid" over the other, which doesn't change the
topology of the diagram.  We will show that this weight system is
the same as the the determinant of the intersection graph of the
diagram.  Our proof is essentially the same as that in \cite{bg};
we include it for completeness.

\begin{thm} \label{T:conway}
\cite{bg} For any chord diagram D, $C(D) = det(\Gamma(D))$.
\end{thm}
{\sc Proof:}  Since both of these weight systems satisfy the
2-term relations, it suffices to show that they agree on caravans.
Consider a caravan $D$ with $m_1$ one-humped camels and $m_2$
two-humped camels, as shown in Figure \ref{F:caravan}.  Then
$adj(\Gamma(D)) \cong [0]^{m_1} \oplus \left[{\matrix{0 & 1 \cr 1
& 0}}\right]^{m_2}$.  So $det(\Gamma(D)) = 0^{m_1}1^{m_2} =
\left\{{\matrix{1\ if\ m_1 = 0 \cr 0\ otherwise}}\right.$.

On the other hand, if we surger all the chords of the diagram, we
obtain an unlink with $m_1 + 1$ components, which means that $C(D)
= \left\{{\matrix{1\ if\ m_1 = 0 \cr 0\ otherwise}}\right.$.  So
the two weight systems agree.  $\Box$

We now turn to the HOMFLYPT polynomial.  We will begin by
considering a framed version of the HOMFLYPT polynomial - i.e. an
invariant of regular isotopy, rather than isotopy.  This invariant
is the Laurent polynomial $P(l,m) \in {\bf Z}[l^{\pm 1},m^{\pm
1}]$ defined by the following skein relations (see \cite{li})
($L^+$ is the result of adding a positive kink to the link $L$):
$$P(L_+)-P(L_-) = mP(L_0)$$
$$P(L^+) = lP(L)$$
$$P(L \cup O) = \frac{l-l^{-1}}{m}P(L)$$
$$P(O) = 1$$
If we make the substitutions $m = e^{ax/2}-e^{-ax/2}$ and $l =
e^{bx/2}$, and expand the resulting power series, we transform the
HOMFLYPT polynomial into a power series in $x$, whose coefficients
are finite type invariants (of regular isotopy).  These invariants
give rise to regular weight systems which we can collect together
as the HOMFLYPT regular weight system $H$.  The skein relations
above give rise to the following relations for $H$, by looking at
the first terms of the power series (as before, $D_v$ is the
result of surgering the chord $v$ in $D$):
$$H(D) = aH(D_v)$$
$$H(D \cup O) = bH(D)$$
$$H(O) = 1$$
So if $D$ is an unlink of $k$ components, $H(D) = b^{k-1}$.  Since
the first of these relations is almost the same as for the Conway
weight system $C$, the same argument shows that $H$ satisfies the
2-term relations.  We will use this to show that the HOMFLYPT
regular weight system is equivalent to the rank of the
intersection graph of the diagram.  (This result was found
independently by Soboleva, for the case $a = 1$ \cite{so}.)

\begin{thm} \label{T:homfly}
For any chord diagram D of degree k, $H(D) =
a^kb^{k-rank(\Gamma(D))} = (ab)^kR(\Gamma(D))(b^{-1})$
\end{thm}
{\sc Proof:}  As with Theorem \ref{T:conway}, it suffices to show
that the weight systems agree on caravans.  Let $D$ be the caravan
with $m_1$ one-humped camels and $m_2$ two-humped camels, as in
Figure \ref{F:caravan} (so the degree of $D$ is $m_1+2m_2$).  As
before, $adj(\Gamma(D)) \cong [0]^{m_1} \oplus \left[{\matrix{0 &
1 \cr 1 & 0}}\right]^{m_2}$, so the rank is $2m_2$.  On the other
hand, if we surger all the chords (each time multiplying $H$ by
$a$), the resulting link has $m_1+1$ components, so $H(D) =
a^kb^{m_1} = a^kb^{k-rank(\Gamma(D))}$. $\Box$

\begin{cor} \label{C:rank}
If D is a chord diagram of degree k, and $L_D$ is the link with c
components obtained by surgering all of the chords of D, then
$rank(\Gamma(D)) = k-c+1$.
\end{cor}

We can also consider the unframed HOMFLYPT polynomial
$\widehat{P}(l,m)$, defined by $\widehat{P}(L) =
l^{-writhe(L)}P(L)$ (see \cite{li}).  This invariant is determined
by the skein relations:
$$l\widehat{P}(L_+)-l^{-1}\widehat{P}(L_-) = m\widehat{P}(L_0)$$
$$\widehat{P}(L \cup O) = \frac{l-l^{-1}}{m}\widehat{P}(L)$$
$$\widehat{P}(O) = 1$$
After making the same substitutions as before, we again obtain a
power series whose coefficients are finite type invariants (this
time of isotopy). The collection of the associated weight systems
$\widehat{H}$ was described by Meng \cite{men} (here $D_v$ is the
result of surgery on the chord $v$, and $D\backslash v$ is the
result of removing the chord $v$):
$$\widehat{H}(D) = a\widehat{H}(D_v) - b\widehat{H}(D\backslash v)$$
$$\widehat{H}(D \cup O) = b\widehat{H}(D)$$
$$\widehat{H}(O) = 1$$

It is easy to see that this weight system is simply the canonical
projection of $H$, and so we can conclude that:

\begin{thm} \label{T:homfly2}
For any chord diagram D of degree k, $\widehat{H}(D) =
(ab)^k\widehat{R}(\Gamma(D))(b^{-1})$
\end{thm}
{\sc Proof:}  Both weight systems are the canonical projections of
$H$. $\Box$

{\it Remark:}  Rather than considering the rank of the adjacency
matrix, we could as easily have studied its nullity.  If we define
$N(G)(x) = x^{nullity(adj(G))}$, and let $\widehat{N}(G)$ be its
canonical projection, then Theorems \ref{T:homfly} and
\ref{T:homfly2} imply that $H(D) = a^kN(\Gamma(D))(b)$ and
$\widehat{H}(D) = a^k\widehat{N}(\Gamma(D))(b)$.

\section{Marked Chord Diagrams and the Kauffman weight system} \label{S:marked}

\begin{figure}
    $$\includegraphics{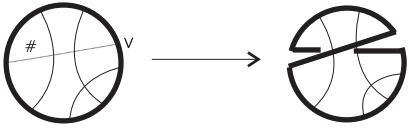}$$
    \caption{Surgery on a marked chord $v$} \label{F:msurgery}
    \end{figure}

In this section we will look at {\it marked} chord diagrams,
motivated by the Kauffman polynomial.  The idea is that, where we
replaced a chord with a band in the previous section, a marked
chord will be replaced by a {\it twisted} band, as in
Figure~\ref{F:msurgery}.
The two different surgeries correspond to the two resolutions of a
crossing, $L_0$ and $L_\infty$, in Figure \ref{F:skein}.

We will begin by defining marked chord diagrams and graphs,
together with a natural map from the space of chord diagrams
(graphs) to the space of marked chord diagrams (graphs).  We will
define an expanded set of 2-term relations on these spaces, and
show that a modification of the adjacency matrix is invariant
under these relations.  We will use this to construct regular
graph weight systems, and show that one of these systems is
equivalent to the regular weight system associated with the
(framed) Kauffman polynomial.

\subsection{Marked Chord Diagrams and Graphs} \label{SS:marked}

A {\it marking} of a chord diagram $D$ (respectively, graph $G$)
is simply a partition of the set of chords $J(D)$ (resp. vertices
$V(G)$) into two disjoint subsets $J_m$ and $J_u$ ($V_m$ and
$V_u$), where $J_m$ ($V_m$) is the set of {\it marked} chords
(vertices), and $J_u$ ($V_u$) is the set of {\it unmarked} chords
(vertices).  We will typically denote a marked chord by labeling
it with a pound sign (\#).

There is a natural map from the vector space of chord diagrams to
the vector space of marked chord diagrams, simply taking a diagram
to the sum (with signs) of all possible ways of marking it.

\begin{defn} \label{D:M}
Consider a chord diagram D, and a subset J of the set of chords of
D. Let $D^J$ denote the marked chord diagram obtained by marking
all the chords in J.  Then we define a map M from the vector space
of chord diagrams to the vector space of marked chord diagrams by:
$M(D) = \sum_{J \subset J(D)}{(-1)^{|J|}D^J}$.
\end{defn}

We can define a similar map (which we will also denote $M$) from
the vector space of graphs to the space of marked graphs.  There
is an obvious lifting of the map $\Gamma$ from a chord diagram to
its intersection graph to a map from a marked chord diagram to its
marked intersection graph, which we will also denote $\Gamma$
(simply mark the vertices corresponding to the marked chords).
Clearly, for any chord diagram $D$, $M(\Gamma(D)) = \Gamma(M(D))$.

\subsection{Extended 2-term relations} \label{SS:ext2term}

We can extend the 2-term relations from section \ref{SS:2term} to
a set of 2-term relations on the space of marked chord diagrams.
Just as the original 2-term relations were motivated by the idea
of replacing chords by bands, the extensions are motivated by the
idea of replacing marked chords by twisted bands.

\begin{figure}
    $$\includegraphics{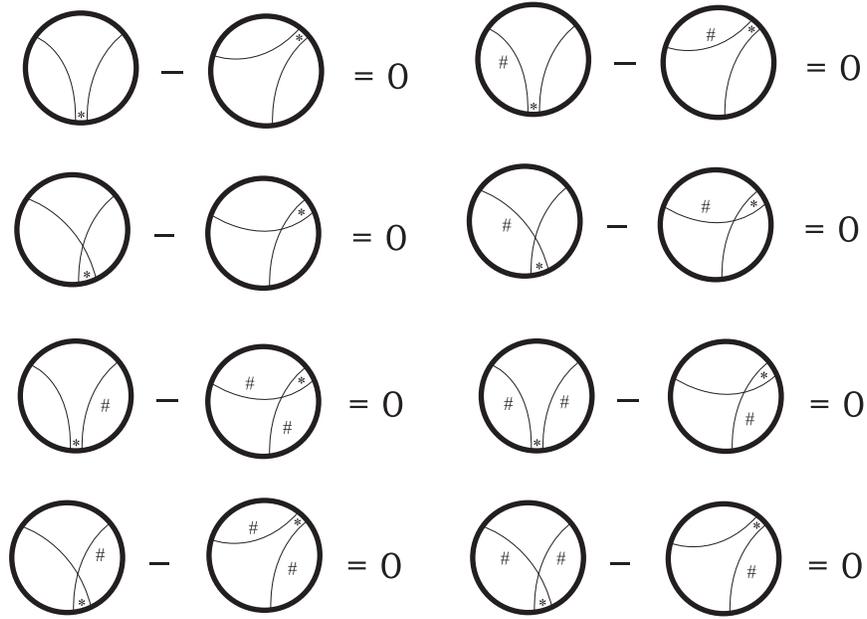}$$
    \caption{2-term relations on marked chord diagrams} \label{F:ext2term}
    \end{figure}

The extended 2-term relations are shown in
Figure~\ref{F:ext2term}.
The space of marked chord diagrams modulo these relations can be
given the structure of a bialgebra by using the same product and
coproduct as in section \ref{SS:finitetype}.  We will denote this
bialgebra (and the underlying graded vector space) $B^m$.  The
only point which needs to be checked is that the product is
well-defined modulo the 2-term relations. This verification is
very similar to the corresponding proof for chord diagrams in
\cite{bn}, and is left as an exercise for the reader.

\begin{prop} \label{P:pullback}
If $\beta$ is a functional on $B^m$, then $\beta \circ M$ is a
regular weight system (a functional on $A$).
\end{prop}
{\sc Proof:}  It is easy to check that the image of a 4-term
relation under $M$ is a linear combination of 2-term relations, so
$M$ is a bialgebra homomorphism from $A$ to $B^m$. $\Box$

We want to find a set of generators for the vector space of marked
chord diagrams modulo the 2-term relations.  One such spanning set
is a generalization of the caravans of the original 2-term
relations.

\begin{defn} \label{D:mcaravan}
A {\bf marked ($n_1, n_2, n_3$)-caravan} is a marked chord diagram
$(\Theta_m)^{n_1}\Theta^{n_2}X^{n_3}$, where $\Theta_m$ is the
chord diagram consisting of a single marked chord (a marked
one-humped-camel), $\Theta$ is the diagram consisting of a single
unmarked chord (a one-humped-camel), and X is the diagram
consisting of two intersecting unmarked chords (a
two-humped-camel).
\end{defn}

\begin{figure}
    $$\includegraphics{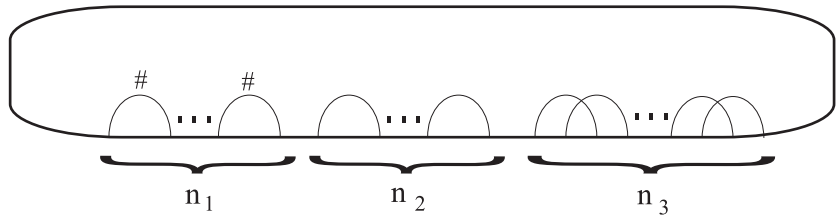}$$
    \caption{A marked ($n_1, n_2, n_3$)-caravan} \label{F:mcaravan}
    \end{figure}

\begin{figure}
    $$\includegraphics{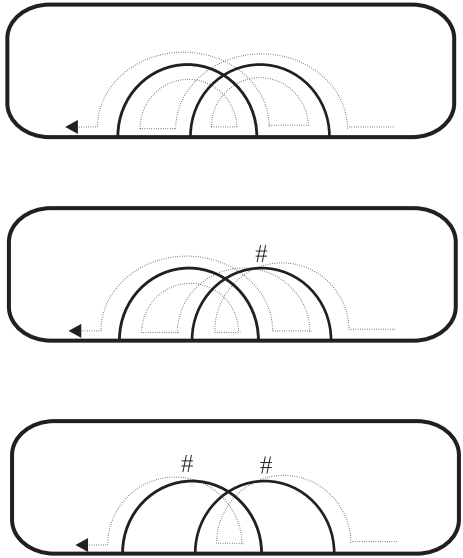}$$
    \caption{Factoring out two-humped-camels.  Notice that a chord
being slid over the camel will follow the same path whether it is
marked or unmarked, and that its marking after being slid over the
camel will be the same as before (although it may change during
the process).} \label{F:2hump}
    \end{figure}

\begin{figure}
    $$\includegraphics{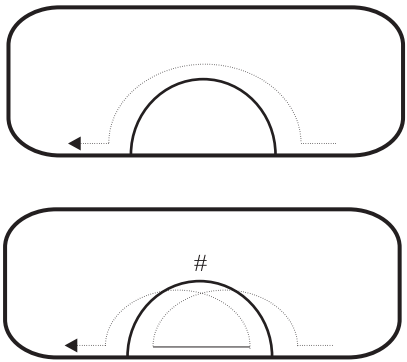}$$
    \caption{Factoring out one-humped-camels.  Notice that a chord
being slid over the camel will follow the same path whether it is
marked or unmarked, and that its marking after being slid over the
camel will be the same as before (although it may change during
the process).} \label{F:1hump}
    \end{figure}

\begin{figure}
    $$\includegraphics{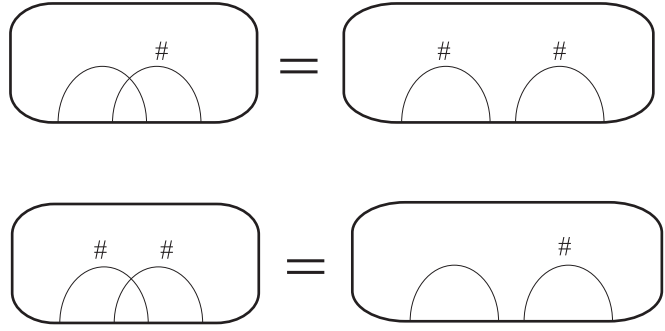}$$
    \caption{Reducing marked two-humped camels to one-humped camels.} \label{F:reduce}
    \end{figure}

An example of a marked caravan is shown in
Figure~\ref{F:mcaravan}.
We will now show that these caravans span the space of marked
chord diagrams, modulo the extended 2-term relations, using an
argument similar to that in \cite{bg}.  In fact, we will show a
slightly stronger fact:

\begin{thm} \label{T:mcaravan}
Any marked chord diagram is equivalent to a marked caravan, modulo
the 2-term relations.
\end{thm}
{\sc Proof:}  Let $D$ be a marked chord diagram.  To begin with,
assume that $D$ has two intersecting chords $c_1$ and $c_2$
(possibly marked). There are four possibilities:  both chords are
unmarked, only $c_1$ is marked, only $c_2$ is marked, or both
chords are marked.  In each case, the pair of chords can be slid
to the right using the 2-term relations, as in
Figure~\ref{F:2hump},
until a (possibly marked) two-humped-camel is factored out.
Continuing inductively, we can factor out (possibly marked)
two-humped camels until there are no remaining pairs of
intersecting chords.  Then, among the remaining chords, there will
be a "smallest" chord, whose endpoints are not separated by the
endpoints of any other chord. This chord, whether marked or
unmarked, can be slid to the right as in Figure~\ref{F:1hump},
until a (possibly marked) one-humped-camel is factored out.
Continuing inductively, we can reduce the remaining chords to a
series of marked and unmarked one-humped-camels.  Finally, we can
reduce the marked two-humped-camels to pairs of one-humped-camels
as in Figure~\ref{F:reduce}.
We are left with a product of marked and unmarked
one-humped-camels and unmarked two-humped-camels, which is a
marked caravan.  This completes the proof. $\Box$

We can define similar 2-term relations for marked graphs. Consider
a marked graph $G$ with vertices $A, B$.  Let $(G)_{AB}$ be this
graph with both $A$ and $B$ unmarked, $(G)_{A*B}$ be the graph
with $A$ marked and $B$ unmarked, $(G)_{AB*}$ be the graph with
$A$ unmarked and $B$ marked, and $(G)_{A*B*}$ be the graph with
both $A$ and $B$ marked.  Then the 2-term relations are ($G'_{AB},
\widetilde{G}_{AB}, \widetilde{G}_{AB}$ are as defined in section
\ref{SS:Lando}):
$$(G)_{AB}-(\widetilde{G}_{AB})_{AB}=0$$
$$(G'_{AB})_{AB}-(\widetilde{G}'_{AB})_{AB}=0$$
$$(G)_{A*B}-(\widetilde{G}_{AB})_{A*B}=0$$
$$(G'_{AB})_{A*B}-(\widetilde{G}'_{AB})_{A*B}=0$$
$$(G)_{AB*}-(\widetilde{G}'_{AB})_{A*B*}=0$$
$$(G'_{AB})_{AB*}-(\widetilde{G}_{AB})_{A*B*}=0$$
$$(\widetilde{G}_{AB})_{AB*}-(G'_{AB})_{A*B*}=0$$
$$(\widetilde{G}'_{AB})_{AB*}-(G)_{A*B*}=0$$
We let $E^m$ denote the vector space of marked graphs modulo these
relations.  Then $E^m$ can be given the structure of a bialgebra
by using the same product and coproduct as in section
\ref{SS:Lando}.

\begin{prop} \label{P:pullback2}
If $\gamma$ is a functional on $E^m$, then $\gamma \circ M$ is a
regular graph weight system, and $\gamma \circ \Gamma$ is a
functional on $B^m$.
\end{prop}
{\sc Proof:}  To show the first part of the proposition, we simply
need to check that the image of a 4-term relation under $M$ is a
linear combination of 2-term relations, so $M$ is a bialgebra
homomorphism from $F$ to $E^m$. The second part of the proposition
is immediate. $\Box$

The commutative diagram below summarizes the maps between the
various bialgebras we have discussed. All of the maps are
bialgebra homomorphisms.  It is worth noting that the map $M$ is
{\it not} a homomorphism from $B$ to $B^m$, because the image of a
2-term relation in $B$ may not be a sum of 2-term relations in
$B^m$. The maps $p$ are the natural projections from $A$ and $F$
to $B$ and $E$, respectively. The maps $\widetilde{p}$ are
projections from $B^m$ and $E^m$ to $B$ and $E$ (respectively),
defined by sending all diagrams (graphs) with marked chords
(vertices) to 0.
$$\includegraphics{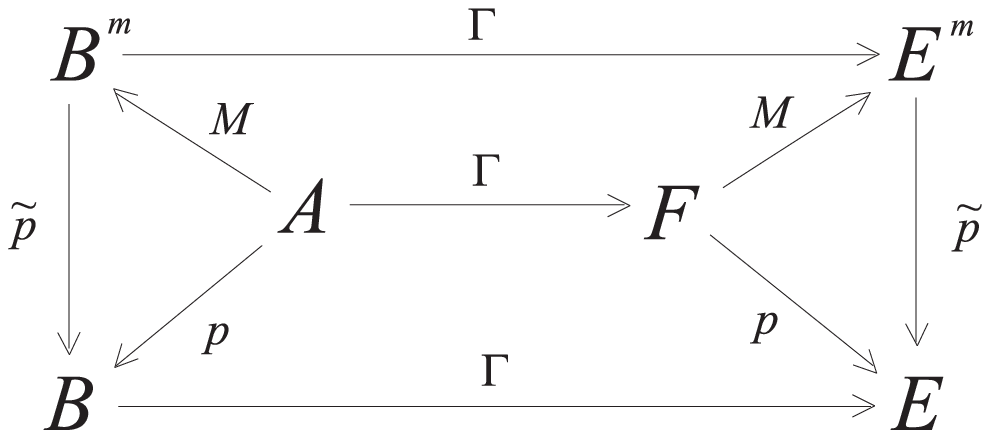}$$

\subsection{Marked Adjacency Matrices} \label{SS:madj}

Now that we have defined the algebra $E^m$, and shown that
functionals on this algebra give rise to regular graph weight
systems and hence (via the deframing map) regular weight systems,
we want to construct explicit examples.  Once again, we will use
the adjacency matrix of a graph.  The adjacency matrix of a marked
graph is defined as in section \ref{SS:adjweight}, except that
$adj(G)_{ii} = 1$ if $v_i$ is a marked vertex.  (We can visualize
a marked vertex as having a small loop attached to it, so it is
adjacent to itself.)

As in Section~\ref{SS:adjweight}, this matrix can be viewed as a
symmetric bilinear form over ${\bf Z}_2$, and is well-defined up
to isomorphism of forms.  As before, we define the $rank$ (resp.
$det$) of a marked graph as the rank (resp. determinant) of the
adjacency matrix of the graph.

\begin{thm} \label{T:madj}
The isomorphism class of the adjacency matrix of a marked graph
satisfies the extended 2-term relations.
\end{thm}
{\sc Proof:}  Consider a graph $G$ with vertices $A$ and $B$.  We
simply need to verify the eight 2-term relations for the adjacency
matrix.  We can do this by writing down the matrices explicitly,
as we did in Theorem \ref{T:weight}. As before, we write our
matrices with respect to the basis $\{A, B, S_{AB}, S_A, S_B,
S_0\}$, and we assume that $A$ and $B$ are connected by an edge in
$G$.  Also as before, $I$ and $O$ represent a row or column of 1's
and 0's respectively:
$$adj((G)_{AB}) = \left[{\matrix{0 & 1 & I & I & O & O \cr
    1 & 0 & I & O & I & O \cr
    I & I & * & * & * & * \cr
    I & O & * & * & * & * \cr
    O & I & * & * & * & * \cr
    O & O & * & * & * & *}}\right] \cong \left[{\matrix{0 & 1 & O & I & I & O \cr
    1 & 0 & I & O & I & O \cr
    O & I & * & * & * & * \cr
    I & O & * & * & * & * \cr
    I & I & * & * & * & * \cr
    O & O & * & * & * & *}}\right] = adj((\widetilde{G}_{AB})_{AB})$$
$$adj((G'_{AB})_{AB}) = \left[{\matrix{0 & 0 & I & I & O & O \cr
    0 & 0 & I & O & I & O \cr
    I & I & * & * & * & * \cr
    I & O & * & * & * & * \cr
    O & I & * & * & * & * \cr
    O & O & * & * & * & *}}\right] \cong \left[{\matrix{0 & 0 & O & I & I & O \cr
    0 & 0 & I & O & I & O \cr
    O & I & * & * & * & * \cr
    I & O & * & * & * & * \cr
    I & I & * & * & * & * \cr
    O & O & * & * & * & *}}\right] = adj((\widetilde{G}'_{AB})_{AB})$$
$$adj((G)_{A*B}) = \left[{\matrix{1 & 1 & I & I & O & O \cr
    1 & 0 & I & O & I & O \cr
    I & I & * & * & * & * \cr
    I & O & * & * & * & * \cr
    O & I & * & * & * & * \cr
    O & O & * & * & * & *}}\right]
\cong \left[{\matrix{1 & 1 & O & I & I & O \cr
    1 & 0 & I & O & I & O \cr
    O & I & * & * & * & * \cr
    I & O & * & * & * & * \cr
    I & I & * & * & * & * \cr
    O & O & * & * & * & *}}\right] = adj((\widetilde{G}_{AB})_{A*B})$$
$$adj((G'_{AB})_{A*B}) = \left[{\matrix{1 & 0 & I & I & O & O \cr
    0 & 0 & I & O & I & O \cr
    I & I & * & * & * & * \cr
    I & O & * & * & * & * \cr
    O & I & * & * & * & * \cr
    O & O & * & * & * & *}}\right]
\cong \left[{\matrix{1 & 0 & O & I & I & O \cr
    0 & 0 & I & O & I & O \cr
    O & I & * & * & * & * \cr
    I & O & * & * & * & * \cr
    I & I & * & * & * & * \cr
    O & O & * & * & * & *}}\right] = adj((\widetilde{G}'_{AB})_{A*B})$$
$$adj((G)_{AB*}) = \left[{\matrix{0 & 1 & I & I & O & O \cr
    1 & 1 & I & O & I & O \cr
    I & I & * & * & * & * \cr
    I & O & * & * & * & * \cr
    O & I & * & * & * & * \cr
    O & O & * & * & * & *}}\right]
\cong \left[{\matrix{1 & 0 & O & I & I & O \cr
    0 & 1 & I & O & I & O \cr
    O & I & * & * & * & * \cr
    I & O & * & * & * & * \cr
    I & I & * & * & * & * \cr
    O & O & * & * & * & *}}\right] = adj((\widetilde{G}'_{AB})_{A*B*})$$
$$adj((G'_{AB})_{AB*}) = \left[{\matrix{0 & 0 & I & I & O & O \cr
    0 & 1 & I & O & I & O \cr
    I & I & * & * & * & * \cr
    I & O & * & * & * & * \cr
    O & I & * & * & * & * \cr
    O & O & * & * & * & *}}\right]
\cong \left[{\matrix{1 & 1 & O & I & I & O \cr
    1 & 1 & I & O & I & O \cr
    O & I & * & * & * & * \cr
    I & O & * & * & * & * \cr
    I & I & * & * & * & * \cr
    O & O & * & * & * & *}}\right] = adj((\widetilde{G}_{AB})_{A*B*})$$
$$adj((\widetilde{G}'_{AB})_{AB*}) = \left[{\matrix{0 & 0 & O & I & I & O \cr
    0 & 1 & I & O & I & O \cr
    O & I & * & * & * & * \cr
    I & O & * & * & * & * \cr
    I & I & * & * & * & * \cr
    O & O & * & * & * & *}}\right]
\cong \left[{\matrix{1 & 0 & I & I & O & O \cr
    0 & 1 & I & O & I & O \cr
    I & I & * & * & * & * \cr
    I & O & * & * & * & * \cr
    O & I & * & * & * & * \cr
    O & O & * & * & * & *}}\right] = adj((G)_{A*B*})$$
$$adj((\widetilde{G}_{AB})_{AB*}) = \left[{\matrix{0 & 1 & O & I & I & O \cr
    1 & 1 & I & O & I & O \cr
    O & I & * & * & * & * \cr
    I & O & * & * & * & * \cr
    I & I & * & * & * & * \cr
    O & O & * & * & * & *}}\right]
\cong \left[{\matrix{1 & 0 & I & I & O & O \cr
    0 & 1 & I & O & I & O \cr
    I & I & * & * & * & * \cr
    I & O & * & * & * & * \cr
    O & I & * & * & * & * \cr
    O & O & * & * & * & *}}\right] = adj((G'_{AB})_{A*B*})$$
The isomorphisms are just the result of adding the second row (and
column) of the matrix on the left to its first row (and column),
modulo 2. So the isomorphism classes of the adjacency matrices of
the graphs satisfy the extended 2-term relations. $\Box$

\begin{cor} \label{C:mrank}
The rank and determinant of a marked graph are functionals on
$E^m$.
\end{cor}

We can combine these functionals with $M$ to obtain regular graph
weight systems.  In order to construct polynomial-valued weight
systems, we will begin with $s(G)=x^{rank(G)}$ and
$t(G)=x^{det(G)}$, whose linear extensions are also functionals on
$E^m$:

\begin{thm} \label{T:s(g)}
Given an unmarked graph G, and a subset J $\subset$ V(G), define
$G^J$ as the result of marking the vertices in J.  Define the maps
S(G) and T(G) as follows:
$$S(G)(x) = \sum_{J \subset V(G)}{(-1)^{|J|}x^{rank(G^J)}}$$
$$T(G)(x) = \sum_{J \subset V(G)}{(-1)^{|J|}x^{det(G^J)}}$$
Then these maps are regular graph weight systems.  Moreover, the
map S(G) is multiplicative:  $S(G_1 \cdot G_2) = S(G_1)S(G_2)$
\end{thm}
{\sc Proof:}  $S(G) = s(M(G))$ and $T(G) = t(M(G))$ (where $s$ and
$t$ are extended linearly), so these are regular graph weight
systems by Proposition \ref{P:pullback2}. Since $rank(G_1 \cdot
G_2) = rank(G_1)+rank(G_2)$, we can see that $s(G_1 \cdot G_2) =
s(G_1)s(G_2)$.  It is easy to check that $M$ is also
multiplicative.  Therefore $S(G)$ is multiplicative. $\Box$

Moreover, we can obtain graph weight systems by applying the
canonical projection from section \ref{SS:Lando}.

\begin{thm} \label{T:s^(g)}
Given an unmarked graph G with n vertices, a subset J $\subset$
V(G), and a subset $J_m \subset J$, define $G^{J_m}_J$ as the
subgraph induced by J, with the vertices in $J_m$ marked.  Then we
define $\widehat{S}(G)$ and $\widehat{T}(G)$ as follows:
$$\widehat{S}(G)(x) = \sum_{J \subset V(G)}{(x-1)^{n-|J|}S(G_J)} =
    \sum_{J \subset V(G)}{\sum_{J_m \subset J}{(-1)^{|J_m|}(x-1)^{n-|J|}
    x^{rank(G^{J_m}_J)}}}$$
$$\widehat{T}(G)(x) = \sum_{J \subset V(G)}{T(G_J)} =
    \sum_{J \subset V(G)}{\sum_{J_m \subset J}{(-1)^{|J_m|}x^{det(G^{J_m}_J)}}}$$
These maps are the canonical projections of S(G) and T(G), and so
are graph weight systems.
\end{thm}
{\sc Proof:}  Recall that the deframing map $\phi(G_1 \otimes G_2)
= (-\bullet)^{deg(G_1)} \cdot G_2$.  It is easy to check that the
map $M:F \rightarrow E^m$ is multiplicative, i.e. $M(G_1 \cdot
G_2) = M(G_1)M(G_2)$.  So $M(\phi(G_1 \otimes G_2)) =
M(-\bullet)^{deg(G_1)}M(G_2)$.  Since $s(G)$ is also
multiplicative, we have:
\begin{eqnarray*}
    S(\phi(G_1 \otimes G_2)) & = & s(M(\phi(G_1 \otimes G_2))) \\
                             & = & s(M(-\bullet))^{deg(G_1)}s(M(G_2))\\
                             & = & (x-1)^{deg(G_1)}s(M(G_2))\\
                             & = & (x-1)^{deg(G_1)}S(G_2)
\end{eqnarray*}
From this, it is straightforward to see that the projection of
$S(G)$ is $\sum_{J \subset V(G)}{(x-1)^{n-|J|}S(G_J)}$, as
desired.

On the other hand, the determinant of a graph with any isolated
unmarked chords is 0.  So for $T$ we have (denoting the graph
consisting of a single marked vertex by $\bullet\#$):
\begin{eqnarray*}
    T(\phi(G_1 \otimes G_2)) & = & t(M(\phi(G_1 \otimes G_2)))\\
    & = & t(M(-\bullet)^{deg(G_1)}M(G_2))\\
    & = & t\left({(\bullet\#-\bullet)^{deg(G_1)}\sum_{J \subset
        V(G_2)}{(-1)^{|J|}G_2^J}}\right)\\
    & = & t\left({\sum_J{\sum_{k=0}^{deg(G_1)}{{{deg(G_1)} \choose k}
        (\bullet\#)^k(-\bullet)^{deg(G_1)-k}(-1)^{|J|}G_2^J}}}\right)\\
    & = & \sum_J{\sum_{k=0}^{deg(G_1)}{{{deg(G_1)} \choose k}
        (-1)^{|J|}(-1)^{deg(G_1)-k}t((\bullet\#)^k(\bullet)^{deg(G_1)-k}G_2^J)}}
\end{eqnarray*}
Since $det((\bullet\#)^k(\bullet)^{deg(G_1)-k}G_x^J) =
\left\{{\matrix{det(G_2^J)\ if\ k = deg(G_1) \cr 0\
otherwise}}\right.$, we know that
$t((\bullet\#)^k(\bullet)^{deg(G_1)-k}G_x^J) =
\left\{{\matrix{t(G_2^J)\ if\ k = deg(G_1) \cr 1\
otherwise}}\right.$.  So our equation reduces to:
\begin{eqnarray*}
    T(\phi(G_1 \otimes G_2)) & = & \sum_J{(-1)^{|J|}\left({t(G_2^J)+\sum_{k=0}^{deg(G_1)-1}
        {{{deg(G_1)} \choose k}(-1)^{deg(G_1)-k}}}\right)}\\
    & = & \sum_J{(-1)^{|J|}\left({t(G_2^J)+(1-1)^{deg(G_1)}-1}\right)}\\
    & = & \sum_J{(-1)^{|J|}t(G_2^J)}-\sum_J{(-1)^{|J|}}\\
    & = & t(M(G_2)) - 0 = T(G_2)
\end{eqnarray*}
From this, we can conclude that the projection of $T(G)$ is
$\widehat{T}(G)(x) = \sum_{J \subset V(G)}{T(G_J)}$, as desired.
$\Box$

\subsection{The Kauffman weight system}

We want to show that $S(\Gamma(D))$ and $\widehat{S}(\Gamma(D))$
are the weight systems associated with the Kauffman polynomial. We
will begin by considering a framed version of the Kauffman
polynomial $F(y,z)$, defined by the following skein relations
($L_+$, $L_-$, $L_0$ and $L_\infty$ are as shown in Figure
\ref{F:skein}, and $L^+$ is the result of adding a postive kink to
$L$):
$$F(L_+)-F(L_-) = z(F(L_0)-F(L_\infty)$$
$$F(L^+) = yF(L)$$
$$F(L \cup O) = \left({\frac{y-y^{-1}}{z}+1}\right)F(L)$$
$$F(O) = 1$$
To derive finite type invariants, we make the substitutions $z =
e^{ax/2}-e^{-ax/2}$ and $y = e^{(b-1)x/2}$.  If we then expand the
polynomial as a power series in $x$, the coefficients will be
finite type invariants.  The regular weight system associated with
this collection of invariants is defined by the skein relations
below. Here $D$ is an unmarked chord diagram, $v$ is a chord in
$D$, $D_v$ is the result of replacing $v$ by an untwisted band,
and $D^v$ is the result of replacing $v$ by a band with a
half-twist:
$$K(D) = a(K(D_v) - K(D^v))$$
$$K(D \cup O) = bK(D)$$
$$K(O) = 1$$
Note that, if $D$ is an unlink of $k$ components, then $K(D) =
b^{k-1}$.

Our first task is to show that this regular weight system factors
through the algebra $B^m$. We define a map $K^m:B^m \rightarrow
{\bf Z}[a,b]$ recursively by the following relations, where $D$ is
a marked chord diagram, and $v$ is a chord in $D$:
$$K^m(D) = \left\{{\matrix{aK^m(D_v)\ if\ v\ is\ unmarked \cr
aK^m(D^v)\ if\ v\ is\ marked}}\right.$$\
$$K^m(D \cup O) = bK^m(D)$$
$$K^m(O) = 1$$
Note that, if $D$ is a diagram with no chords and $k$ components,
then $K^m(D) = b^{k-1}$.

\begin{prop} \label{P:K^m}
$K^m$ satisfies the extended 2-term relations.
\end{prop}
{\sc Proof:}  In each of the 2-term relations of Figure
\ref{F:ext2term}, replace each unmarked chord by an untwisted band
and each marked chord by a band with a half-twist.  It is clear
that the relations are simply the result of sliding one band over
another, and don't change the topology of the diagram.  We need
only keep in mind that when a band is slid over a half-twisted
band (marked chord), it receives a half-twist itself.  We view a
band with a full twist as equivalent to an untwisted band, since
it does not change the number of components of the diagram, which
is all that matters in the base case of the definition of $K^m$.
$\Box$

\begin{prop} \label{P:lift}
$K = K^m \circ M$, so K is the pullback of $K^m$ by M.
\end{prop}
{\sc Proof:}  Consider a diagram $D$ in $A$.  We will prove the
proposition via induction on the number of chords of $D$.  If $D$
has no chords, then $M(D) = D$.  Since $K$ and $K^m$ differ only
in their first skein relation (which only applies if there are
chords), we conclude that $K^m(M(D)) = K^m(D) = K(D)$.

For our inductive step, assume $D$ has a chord $v$.  Note that
$D_v$ and $D^v$ each have fewer chords than $D$, so $K^m(M(D_v)) =
K(D_v)$ and $K^m(M(D^v)) = K(D^v)$.  If $J$ is a subset of the
chords of $D$, we let $D^J$ denote the marked chord diagram which
results by marking all the chords in $J$. (However, for the single
chord $v$, we will still let $D^v$ denote the result of replacing
$v$ with a half-twisted band.) $M(D)$ is then given by:
$$M(D) = \sum_J{(-1)^{|J|}D^J} = \sum_{J\ s.t.\ v \notin
J}{(-1)^{|J|}(D^J - D^{J \cup v})}$$

Then:
\begin{eqnarray*}
    K^m(M(D)) & = & \sum_{J\ s.t.\ v \notin J}{(-1)^{|J|}(K^m(D^J) - K^m(D^{J \cup
v}))}\\
    & = & \sum_{J\ s.t.\ v \notin J}{(-1)^{|J|}(aK^m((D_v)^J) -
    aK^m((D^v)^J))}\\
    & = & a(K^m(M(D_v)) - K^m(M(D^v)))\\
    & = & a(K(D_v) - K(D^v))\\
    & = & K(D)
\end{eqnarray*}
So by induction, we conclude that for any diagram $D$, $K(D) =
K^m(M(D))$.  $\Box$

\begin{thm} \label{T:kauff}
For any $D \in A$ of degree $k$, $K(D) =
(ab)^kS(\Gamma(D))(b^{-1})$.
\end{thm}
{\sc Proof:}  Since $S(\Gamma(D)) = s(M(\Gamma(D)) =
s(\Gamma(M(D))$, and $K(D) = K^m(M(D))$, it suffices to show that
$(ab)^k(s \circ \Gamma(D))(b^{-1}) = K^m(D)$ for any $D \in B^m$.
Since both of these maps satisfy the extended 2-term relations, it
suffices to show that they agree on marked caravans, by Theorem
\ref{T:mcaravan}.

Consider a marked ($n_1, n_2, n_3$)-caravan $D$, as shown in
Figure \ref{F:mcaravan}. The degree of this caravan is $k = n_1 +
n_2 + 2n_3$. Then $adj(\Gamma(D)) \cong [1]^{n_1} \oplus [0]^{n_2}
\oplus \left[{\matrix{0 & 1 \cr 1 & 0}}\right]^{n_3}$. So
$rank(\Gamma(D)) = n_1 + 2n_3$, and $(ab)^ks(\Gamma(D))(b^{-1}) =
(ab)^kb^{-n_1-2n_3} = a^kb^{k-n_1-2n_3} = a^kb^{n_2}$.

On the other hand, $K^m(D)$ is computed by replacing all the
unmarked chords with untwisted bands and all the marked chords
with twisted bands (multiplying by $a$ each time), and then
looking at the number of components of the resulting link.  This
link will have $n_2+1$ components, so $K^m(D) = a^kb^{n_2} =
(ab)^ks(\Gamma(D))(b^{-1})$, which completes the proof. $\Box$

We can also consider the unframed Kauffman polynomial
$\widehat{F}(y,z)$, defined by $\widehat{F}(L) =
y^{-writhe(L)}F(L)$ (see \cite{li}). This invariant is also
determined by the skein relations:
$$y\widehat{F}(L_+)-y^{-1}\widehat{F}(L_-) = m(\widehat{F}(L_0)-\widehat{F}(L_\infty))$$
$$\widehat{F}(L \cup O) = \left({\frac{y-y^{-1}}{z}+1}\right)\widehat{F}(L)$$
$$\widehat{F}(O) = 1$$
After making the same substitutions as before, we again obtain a
power series whose coefficients are finite type invariants (this
time of isotopy). The collection of the associated weight systems
$\widehat{K}$ was described by Meng \cite{men} (here $D_v$ is the
result of replacing the chord $v$ by an untwisted band, $D^v$ is
the result of replacing the chord $v$ by a half-twisted band, and
$D\backslash v$ is the result of removing the chord $v$):
$$\widehat{K}(D) = a\widehat{K}(D_v) - a\widehat{K}(D^v) - b\widehat{K}(D\backslash v)$$
$$\widehat{K}(D \cup O) = b\widehat{K}(D)$$
$$\widehat{K}(O) = 1$$

It is easy to see that this weight system is simply the canonical
projection of $K$, and so we can conclude that:

\begin{thm} \label{T:kauff2}
For any chord diagram D of degree k, $\widehat{K}(D) =
(ab)^k\widehat{S}(\Gamma(D))(b^{-1})$
\end{thm}
{\sc Proof:}  Both weight systems are the canonical projections of
$K$. $\Box$

{\it Remark:}  Rather than considering the rank of the marked
adjacency matrix, we could as easily have studied its nullity. If
we define $u(G)(x) = x^{nullity(adj(G))}$ and $U(G) = u(M(G))$,
and let $\widehat{U}(G)$ be the canonical projection of $U(G)$,
then Theorems \ref{T:kauff} and \ref{T:kauff2} imply that $K(D) =
a^kU(\Gamma(D))(b)$ and $\widehat{K}(D) =
a^k\widehat{U}(\Gamma(D))(b)$.

We now have explicit formulas for computing the Conway, HOMFLYPT
and Kauffman weight systems directly from intersection graphs.
Hopefully, these interpretations will help shed some light on the
geometric meanings of these polynomials.

\section{Acknowledgements}

I would like to thank Dror Bar-Natan and Louis Zulli for informing
me of their previous work.  I would also like to thank Sergei
Lando for sending me E. Soboleva's paper, along with his own work
\cite{la} on obtaining Vassiliev invariants from intersection
graphs.  Finally, I would like to thank the anonymous reviewer who
made several suggestions which led to substantial revisions of the
paper.

\end{document}